\documentclass[10pt, a4paper]{article}
\usepackage[ansinew]{inputenc}
\usepackage[T1]{fontenc}
\usepackage{longtable}
\usepackage{xcolor}
\usepackage{mathrsfs}
\usepackage{amssymb,amsmath,amsthm}
\usepackage[matrix, arrow, curve]{xy}
\usepackage{tikz-cd}
\usepackage{mathtools}
\usetikzlibrary{quotes,babel, angles}

\newcommand{\arrow}[1]{\stackrel{#1}{\longrightarrow}}

\begin{document}
\setlength{\parindent}{0em}

\title{A necessary and sufficient condition for the existence of non-trivial $S_n$-invariants in the splitting algebra}
\author{Kevin Schlegel}
\date{}

\maketitle

\begin{abstract} For a monic polynomial $f$ over a commutative, unitary ring $A$ the splitting algebra $A_f$ is the universal $A$-algebra such that $f$ splits in $A_f$. The symmetric group acts on the splitting algebra by permuting the roots of $f$. It is known that if the intersection of the annihilators of the elements $2$ and $D_f$ (where $D_f$ depends on $f$) in $A$ is zero, then the invariants under the group action are exactly equal to $A$. We show that the converse holds. 
\end{abstract}

\section*{Introduction}

Given a commutative, unitary ring $A$ and a monic polynomial $f$ of degree $n$ with coefficients in $A$, we can construct the \textbf{splitting algebra} $A_f = A[\tau_1, \tau_2, \dots, \tau_n]$ such that $f$ splits in $A_f$ with the \textbf{universal roots} $\tau_1, \tau_2,\dots, \tau_n$. The symmetric group $S_n$ acts on $A_f$ by permuting the universal roots.\\
One can show that the polynomial ring $A[t_1, t_2,\dots, t_n]$ of degree $n$ over $A$ can be viewed as a splitting algebra with the universal roots $t_1,t_2,\dots, t_n$. Further, if $A$ is a field one can show that the splitting algebra $A_f$ modulo an ideal is isomorphic to the splitting field of $f$. This yields an approach to symmetric functions and Galois theory using splitting algebras (see \cite{EL}).\\
In both cases the elements in $A_f$ that are invariant under the action of $S_n$ are of particular interest. It is always true that the elements in $A\subseteq A_f$ are invariants.\\

\textbf{Question.}\textit{ When exactly are there no elements in $A_f\backslash A$ that are invariant under the action of $S_n$?}\\

As it turns out, this solely depends on the annihilators of the elements $2$ and $D_f:= \prod_{1\leq i < j \leq n} (\tau_i + \tau_j)$ in $A$. As shown in \cite[Theorem 7, p.\ 4]{Th} if their intersection is zero, then the only invariants are $A$. We show that the converse holds as well:\\

\textbf{Main theorem.} \textit{The elements in $A_f$ that are invariant under the action of $S_n$ are exactly the elements in $A$ if and only if}
\begin{align}\tag{$\ast$}
    \text{Ann}_A\,2 \cap \text{Ann}_A\,D_f = 0.
\end{align}

\section{The splitting algebra}

In this section we construct the splitting algebra as in \cite[(1.6), p.\ 30]{PZ}. Further, we follow the conventions in \cite{EL} and restate some results and proofs.\\

\textbf{1.1 Notation and conventions.} We always consider commutative, unitary rings and monic polynomials. Given rings $A$ and $B$, a polynomial
\begin{align*}
    f = t^n + a_1 t^{n-1} + \dots + a_n
\end{align*}
with coefficients in $A$ and a ring homomorphism $\varphi: A \rightarrow B$, we define $\varphi(f)$ as the polynomial
\begin{align*}
    \varphi(f) := t^n + \varphi(a_1)t^{n-1} + \dots + \varphi(a_n)
\end{align*}
with coefficients in $B$.\\
A ring homomorphism $\varphi: A \rightarrow B$ will be called an $A$-algebra. We also call $B$ an $A$-algebra (with the morphism $\varphi$ in mind). Given $A$-algebras $\varphi: A \rightarrow B$ and \linebreak$\psi: A\rightarrow C$, we call a ring homomorphism $\alpha: B \rightarrow C $ an $A$-algebra homomorphism if the diagram
\begin{equation*}
\begin{tikzcd}
A \arrow[d, swap, "\varphi"]  \arrow[dr, "\psi"]&\\
B \arrow[r, swap, "\alpha"] &  C
\end{tikzcd}
\end{equation*}
commutes. We denote the symmetric group of degree $n$ by $S_n$. Further, for $\sigma$ in $S_n$ we sometimes just write $\sigma$ for the group that is generated by $\sigma$.\\

\textbf{1.2 Construction.} Let $A$ be a ring and $f$ a polynomial of degree $n$ with coefficients in $A$. We construct, by descending induction on $i$, $A$-algebras $A_i$ and polynomials $f_i$ with coefficients in $A_i$ for $i=n, n-1, \dots, 0$ as follows:\\
We start with $A_n := A$ and $f_n := f$. Assume that $A_i$ and $f_i$ are already constructed. We define $\tau_i$ as the class of $t$ in $A_i[t]/(f_i)$ and set
\begin{align*}
    \begin{split}
        A_{i-1} &:= A_i[\tau_i]\\
        f_{i-1} &:= f_i/(t-\tau_i).
    \end{split}
\end{align*}
Then $A_f := A_0 = A[\tau_n, \tau_{n-1}, \dots, \tau_1]$ is the \textbf{splitting algebra (of $f$ over $A$)}. By construction $f$ splits in $A_f$ with the \textbf{universal roots (of $f$)} $\tau_1, \tau_2, \dots, \tau_n$.\\

In the following we fix the ring $A$, the polynomials $f_1, f_2, \dots, f_n = f$ and the universal roots $\tau_1, \tau_2, \dots, \tau_n$ as in the above definition.\\

\textbf{1.3 Remark.} Let $\tau_1', \tau_2', \dots, \tau_{n-1}'$ be the universal roots of $f_{n-1}$ in $A[\tau_n]_{f_{n-1}}$. By the construction, there exists a canonical isomorphism of $A$-algebras
\begin{align*}
A[\tau_n]_{f_{n-1}}\longrightarrow A_f    
\end{align*}
via $\tau_{i}' \mapsto \tau_i$ for $i=1, 2, \dots, n-1$. Thus, we may identify $A[\tau_n]_{f_{n-1}} = A_f$ and inductively $A[\tau_{i+1}, \tau_{i+2}, \dots, \tau_{n}]_{f_i} = A_f$ for $i> 0$.\\

\textbf{1.4 Remark.} By the construction, $A_f$ is a free $A$-module with the basis 
\begin{align*}
    B:= \{\tau_1^{m_1}\tau_2^{m_2} \dots \tau_n^{m_n}\mid 0\leq m_i < i\text{ for }1\leq i\leq n\}.
\end{align*}
We call $B$ the \textbf{standard basis (of $A_f$)}. Notice, we may consider $A_f$ as a free $A[\tau_{i+1},\tau_{i+2},\dots, \tau_n]$-module for $i>0$ and then the standard basis of $A[\tau_{i+1},\tau_{i+2}, \dots, \tau_n]_{f_i} = A_f$ is different for all $i$.\\

\textbf{1.5 The universal property.} \textit{For every $A$-algebra $\varphi:A\rightarrow B$ such that $\varphi(f)$ splits in $B$ with some roots $\nu_1, \nu_2, \dots, \nu_n$, there is an $A$-algebra homomorphism $\psi: A_f \rightarrow B$ that is uniquely determined by $\tau_i \mapsto \nu_i$ for $i=1,2,\dots, n$.}\\

\textbf{Proof.} We prove the existence of $\psi$ by induction on $n$. For $n = 1$ we obtain $A_f = A$ and chose $\psi := \varphi$.\\
For $n>1$ consider the $A$-algebra homomorphism $\varphi': A[\tau_n] \rightarrow B$ defined by $\tau_n\mapsto \nu_n$. Since $\varphi(f)(\nu_n) = 0$, $\varphi'$ is well-defined. Hence, we can view $B$ as an $A[\tau_n]$-algebra and $\varphi'(f_{n-1})$ splits in $B$ with the roots $\nu_1, \nu_2, \dots, \nu_{n-1}$. By induction, there exists an $A[\tau_n]$-algebra homomorphism
\begin{align*}
    \psi: A[\tau_n]_{f_{n-1}} \longrightarrow B
\end{align*}
determined by $\tau_i \mapsto \nu_{i}$ for $i=1,2, \dots, n-1$. Since $A[\tau_n]_{f_{n-1}}=A_f$, it follows that $\psi$ is the desired $A$-algebra homomorphism. Because $\tau_1, \tau_2, \dots, \tau_n$ generate $A_f$ as an $A$-algebra, $\psi$ is unique. \qed\\

\textbf{1.6 Definition.} For $\sigma$ in $S_n$ the polynomial $f$ splits in $A_f$ with the roots $\tau_{\sigma(1)}, \tau_{\sigma(2)}, \dots, \tau_{\sigma(n)}$. By the universal property of $A_f$, we obtain an $A$-algebra homomorphism $\psi_{\sigma} : A_f \rightarrow A_f$ determined by $\tau_i \mapsto \tau_{\sigma(i)}$ for $i=1,2,\dots, n$.\linebreak This defines a (left) group action of $S_n$ on $A_f$ by $\sigma x := \psi_\sigma(x)$ for $\sigma $ in $S_n$ and $x$ in $A_f$.\\
Let $G$ be a subgroup of $S_n$. We define the $A$-Algebra $A_f^G$ as the set of all elements $x$ in $A_f$ with $\sigma x = x$ for all $\sigma$ in $G$. Further, the elements in $A_f^G$ are called \textbf{$G$-invariants}.\\

\textbf{1.7 Remark.} Let $G$ be a subgroup of $S_n$. For $\sigma$ in $S_n$ we obtain
\begin{align*}
    \sigma A_f^G = A_f^{\sigma G \sigma^{-1}}.
\end{align*}

\textbf{1.8 Remark.} Let $G \subseteq H$ be subgroups of $S_n$. Then $A_f^{H}\subseteq A_f^{G}$. Further, we can lift $G$-invariants to $H$-invariants via the $A$-module homomorphism
\begin{align*}
    \Sigma_{G}^H : A_f^G \longrightarrow A_f^H, x\mapsto \sum_{\sigma G \in H/G} \sigma x.
\end{align*}

\section{The generic splitting algebra}

In this section the polynomial ring $A[t_1, t_2,\dots, t_n]$ of degree $n$ is constructed as a splitting algebra as in $\cite{EL}$. By that, we are able to introduce a new concept of \textbf{trivial $G$-invariants} for an arbitrary subgroup $G$ of $S_n$ that will be important for the proof of the main theorem.\\

\textbf{2.1 Definition.} Let $g_1, g_2, \dots, g_n$ be independent variables over $A$ and consider the polynomial $g:=t^n + g_1 t^{n-1} + \dots + g_n$. The \textbf{generic splitting algebra (of degree $n$ over $A$)} is the splitting algebra of $g$ over $A[g_1, g_2, \dots, g_n]$ and it is denoted by $\Omega_n(A)$.\\
Further, we denote the universal roots of $\Omega_n(A)$ by $\omega_1, \omega_2, \dots, \omega_n$ and for a subgroup $G$ of $S_n$ we write $\Omega_n^G(A)$ for the $G$-invariants in $\Omega_n(A)$.\\

\textbf{2.2 Lemma.} \textit{Let $B$ and $C$ be $A$-algebras, $\varphi: B\rightarrow C$ an $A$-algebra homomorphism and $h$ a polynomial of degree $n$ with coefficients in $B$. If $\varphi(h)$ splits with some roots $\nu_1, \nu_2, \dots, \nu_n$, then there is a commutative diagram of $A$-algebra homomorphisms}
\begin{equation*}
    \begin{tikzcd}
{A[g_1, g_2, \dots, g_n]} \arrow[r, "\alpha"] \arrow[d, hook] & B \arrow[d, "\varphi"] \\
\Omega_n(A) \arrow[r, "\psi"']                              & C                     
\end{tikzcd}
\end{equation*}
\textit{that is uniquely determined by $\alpha(g) = h$ and $\psi(\omega_i) = \nu_i$ for all $i$.}\\

\textbf{Proof.} Because $g_1, g_2, \dots, g_n$ are independent variables over $A$, we may chose $\alpha$ such that $\alpha(g) =h$. Now consider $C$ as an $A[g_1, g_2, \dots, g_n]$-algebra with the morphism $\varphi \alpha$. Since $\varphi \alpha (g) = \varphi(h)$ splits with the roots $\nu_1, \nu_2, \dots, \nu_n,$ we obtain the desired $\psi$ by the universal property of splitting algebras. \qed\\

\textbf{2.3 Proposition.} \textit{Let $A[t_1, t_2, \dots, t_n]$ be the polynomial ring of degree $n$ over $A$.\linebreak Then there is an $A$-algebra isomorphism}
\begin{align*}
    \beta: A[t_1, t_2, \dots, t_n] \longrightarrow \Omega_n(A), t_i \mapsto \omega_i. 
\end{align*}

\textbf{Proof.} We chose $B:=C:= A[t_1, t_2, \dots, t_n]$ and $h := (t-t_1)(t-t_2)\dots (t-t_n)$
in Lemma 2.2. This yields an $A$-algebra homomorphism
\begin{align*}
    \psi: \Omega_n(A) \longrightarrow A[t_1, t_2, \dots, t_n], \omega_i \mapsto t_i.
\end{align*}
Clearly, $\beta$ and $\psi$ are mutually inverse.\qed\\

\textbf{2.4 Definition.} We choose $B:= A, C:= A_f$ and $h:=f$ in Lemma 2.2 to obtain an $A$-algebra homomorphism 
\begin{align*}
\psi: \Omega_n(A) \longrightarrow A_f, \omega_i \mapsto \tau_i.    
\end{align*}
For a subgroup $G$ of $S_n$ we define the $A$-algebra $A_{f}^{G, \text{tr}}:= \psi \Omega_n^G$ and call the elements in $A_f^{G,\text{tr}}$ \textbf{trivial $G$-invariants}. Clearly, $\psi$ and $\sigma$ commute for all $\sigma$ in $S_n$ and it follows that $A_f^{G,\text{tr}}$ is contained in $A_f^G$. We call the elements in $A_f^G\backslash A_f^{G,\text{tr}}$ \textbf{non-trivial $G$-invariants.}\\

In the following we fix $\psi:\Omega_n(A) \longrightarrow A_f$ as in the above definition.\\

\textbf{2.5 Remark.} Let $G\subseteq H$ be subgroups of $S_n$ and let $\sigma$ be in $S_n$. Then
\begin{align*}
\begin{split}
\sigma A_f^{G, \text{tr}} &= A_f^{\sigma G \sigma^{-1}, \text{tr}}\\
A_f^{H,\text{tr}} &\subseteq A_f^{G,\text{tr}}\\
\Sigma_G^H A_f^{G,\text{tr}} &\subseteq A_f^{H,\text{tr}}.
\end{split}
\end{align*}
\textbf{Proof.} Because $\psi$ and $\sigma$ commute for all $\sigma$ in $S_n$, the claim follows. \qed\\

\textbf{2.6 Remark.} Let $m\leq n$ and $G$ a subgroup of $S_m$. Then
\begin{align*}
    A_f^{G,\text{tr}} = A[\tau_{m+1}, \tau_{m+2}, \dots, \tau_n]_{f_m}^{G,\text{tr}}.
\end{align*}

\textbf{Proof.} Throughout we make use of Proposition 2.3. We have a commutative diagram of $A$-algebra homomorphisms
\begin{equation*}
\begin{tikzcd}
{\Omega_n(A)} \arrow[r, "\beta"] \arrow[d] & \Omega_m(A[\tau_{m+1},\tau_{m+2}, \dots, \tau_n]) \arrow[d] \\
A_f \arrow[r, equal]                              & A[\tau_{m+1},\tau_{m+2}, \dots, \tau_n]_{f_m}                     
\end{tikzcd}
\end{equation*}
with $\beta(\omega_i):= \omega_i$ for $1\leq i\leq m$ and $\beta(\omega_j) := \tau_j$ for $j>m$. We will show
\begin{align*}
    \beta \Omega_n^G(A) = \Omega_m^G(A[\tau_{m+1}, \tau_{m+2}, \dots, \tau_n])
\end{align*}
so the claim follows. Since $\beta$ and $\sigma$ commute for $\sigma$ in $G$ we obtain $"\subseteq"$. Now let $x\in \Omega_m^G(A[\tau_{m+1}, \tau_{m+2}, \dots, \tau_n])$ and write $x=\sum a_p p$, where we sum over all monic monomials $p$ in the variables $\omega_1, \omega_2, \dots, \omega_m$ with coefficients $a_p$ in $A[\tau_{m+1}, \tau_{m+2},\dots, \tau_n]$. Then $\sigma x = x$ for $\sigma$ in $G$ if and only if $a_p = a_{\sigma p}$ for all $p$. For every $a_p$ we can choose $b_p$ in $A[\omega_{m+1}, \omega_{m+2}, \dots, \omega_{n}]$ such that $\beta(b_p) = a_p$ and $b_p = b_{\sigma p}$ for all $\sigma$ in $G$. Let $y:= \sum b_p p$. Then $\beta(y) = x$ and $\sigma y = y$ for all $\sigma \in G$. Hence $"\subseteq"$ follows. \qed

\section{($\ast$) is sufficient}

In this section our goal is to show the sufficiency of $(\ast)$ in the main theorem. It first has been proven in \cite[(3.6), p.\ 49-50]{PZ} in a slightly different formulation (namely $D_f$ is replaced with the discriminant of $f$, however this formulation is equivalent, see (i) $\Leftrightarrow$ (ii) in \cite[Theorem 7, p.\ 4]{Th}). It was also proven in \cite{Th}; we make use of similar methods and combine them with the concept of trivial invariants.\\

\textbf{3.1 Lemma.} \textit{Let $n=2$. Then $\tau_1+\tau_2 \in A$ and $A_f^{S_2}$ decomposes as an $A$-module via}
\begin{align*}
    A_f^{S_2} = A\oplus \left(\text{Ann}_A\,2\cap \text{Ann}_A\, (\tau_1+\tau_2) \right) \tau_2.
\end{align*}
\textit{In particular} $\text{Ann}_{A_f}\,2\cap \text{Ann}_{A_f}\, (\tau_1+\tau_2)\subseteq A_f^{S_2}$.\\

\textbf{Proof.} We write $f = t^2+p t +q$ for $p$ and $q$ in $A$. Since $f$ factors in $A_f$ via $f = (t-\tau_1)(t-\tau_2)$, we obtain $\tau_1+\tau_2 = -p \in A$. Using the standard basis of $A_f$,\linebreak an arbitrary element $x$ in $A_f$ can be written as $x = a+b\tau_2$ with $a$ and $b$ in $A$. Then, $x$ is an $S_2$-invariants if and only if 
\begin{align*}
    0 = x-(1\,2)x = b(\tau_2-\tau_1) = bp+2b\tau_2.
\end{align*}
Now the right hand side equals zero if and only if $b\in \text{Ann}_A\,p\cap \text{Ann}_A\,2$. This proves the claim. \qed \\

\textbf{3.2 Lemma.} \textit{Let $n>2$. Then}
\begin{align*}
    A_f^{S_n} \cap A[\tau_n] = A.
\end{align*}
\textbf{Proof.} Let $p$ be the coefficient of $t^{n-1}$ in $f$ and let $q$ be the coefficient of $t^{n-2}$ in $f_{n-1}$. Since $f = (t-\tau_n)f_{n-1}$ we obtain $q = p +\tau_n$. Further, $\tau_{n-1}$ is a root of $f_{n-1}$ which yields
\begin{align*}
    0 = f_{n-1}(\tau_{n-1}) = \tau_{n-1}^{n-1} + (p+\tau_n){\tau_{n-1}^{n-2}} + \sum_{i=0}^{n-3} q_i \tau_{n-1}^i 
\end{align*}
with $q_i$ in $A[\tau_n]$. Hence, we can replace $\tau_n \tau_{n-1}^{n-2}$ with $\tau_{n-1}^{n-1}$ in the standard basis of $A_f$ and obtain a new basis $B'$. Notice, we have $\tau_n \tau_{n-1}^{n-2} \neq \tau_{n}^{j}$ for all $j\leq n-1$ (because $n>2$). Now let $x = \sum_{i=0}^{n-1} a_i \tau_n^i$ be an $S_n$-invariant with $a_i$ in $A$ and $\sigma := (n-1\, n)$. Then
\begin{align*}
    0 = x - \sigma x = \sum_{i=1}^{n-1} a_i \tau_n^i - \sum_{i=1}^{n-1} a_i \tau_{n-1}^i.
\end{align*}
Since the right hand side is expressed in the basis $B'$, we obtain $a_i = 0$ for all $i>0$.\qed\\

\textbf{3.3 Lemma.} \textit{Let $n\geq i \geq 2$. Then}
\begin{align*}
    A_f^{S_n} \cap A[\tau_{i+1}, \tau_{i+2}\dots, \tau_n] = A.
\end{align*}
\textbf{Proof.} We will show the claim by induction on $n$. For $n=2$ it is trivial. Let $n>2$. Then $A_f^{S_n} \subseteq A[\tau_n]_{f_{n-1}}^{S_{n-1}}$ and by induction
\begin{align*}
    A_f^{S_n} \cap A[\tau_{i+1}, \tau_{i+2}, \dots, \tau_n] \subseteq A[\tau_n]_{f_{n-1}}^{S_{n-1}} \cap A[\tau_n][\tau_{i+1}, \tau_{i+2}, \dots, \tau_{n-1}] = A[\tau_n].
\end{align*}
Now the claim follows by Lemma 3.2.\qed\\

\textbf{3.4 Proposition.} \textit{Let $n\geq 2$. Then }
\begin{align*}
    A_f^{S_n} \cap A_f^{S_2, \textnormal{tr}} = A.
\end{align*}

\textbf{Proof.} First, we proof the claim for $n=2$: It is equivalent to $A_f^{S_2, \text{tr}} = A$. Consider the generic splitting algebra $\Omega_2(A) = A[g_1, g_2]_g$. As seen before, we have $\omega_1+\omega_2=-g_1$ and thus $\text{Ann}_{A[g_1, g_2]}(\omega_1+\omega_2) = 0$. By Lemma 3.1 we obtain $\Omega_2^{S_2}(A) = A[g_1, g_2]$ and the commutative diagram
\begin{equation*}
\begin{tikzcd}
{A[g_1, g_2]} \arrow[r] \arrow[d, hook] & A \arrow[d, hook] \\
\Omega_n(A) \arrow[r, "\Psi"']                              & A_f                     
\end{tikzcd}
\end{equation*}
yields the desired result. For $n>2$, we obtain 
\begin{align*}
A_f^{S_2,\text{tr}} = A[\tau_3, \tau_4, \dots, \tau_n]_{f_2}^{S_2, \text{tr}} = A[\tau_3, \tau_4, \dots, \tau_n].
\end{align*}
By Lemma 3.3 this yields the desired result. \qed \\

\textbf{3.5 Corollary.} \textit{The trivial $S_n$-invariants are exactly equal to $A$.}\\

\textbf{3.6 Definition.} We define $D_f:= \prod_{1\leq i < j \leq n}(\tau_i+\tau_j)$. Since $D_f$ is the image of $\prod_{1\leq i < j \leq n}(\omega_i+\omega_j)\in \Omega_n^{S_n}(A)$ under $\psi$, it follows that $D_f$ is a trivial $S_n$-invariant so $D_f\in A$ by Corollary 3.5.\\

The following proposition proves the sufficiency of $(\ast)$ in the main theorem.\\

\textbf{3.7 Proposition.} \textit{Let $x$ be an $S_n$-invariant. Then, both $2x$ and $D_f x$ are contained in $A$. In particular if $\textnormal{Ann}_A\,2 \cap \textnormal{Ann}_A\,D_f = 0$, then $A_f^{S_n} = A$.}\\

\textbf{Proof.} We have $2x \in A_f^{S_n}$ and $2x = \Sigma_{1}^{(1\,2)} x \in A_f^{(1\,2), \text{tr}}$. By Proposition 3.4 it follows that $2x \in A$.\\
Now we define
\begin{align*}
E_f := \prod_{\substack{1\leq i < j\leq n \\ (i,j)\neq (1,2)}}(\tau_i+\tau_j).    
\end{align*}
Then $D_f = E_f (\tau_1+\tau_2)$ and $E_f\in A_f^{(1\,2),\text{tr}}$. Further, we obtain
\begin{align*}
(\tau_1+\tau_2) x = \Sigma_1^{(1\,2)}(\tau_1x) \in A_f^{(1\,2),\text{tr}}.    
\end{align*}
It follows that $D_f x = E_f (\tau_1+\tau_2)x\in A_f^{(1\,2),\text{tr}}$ but also $D_f x \in A_f^{S_n}$. Again by Proposition 3.4 we obtain $D_f x \in A$. Now if $x\in A_f^{S_n}\backslash A$, then we find a non-zero element in $\text{Ann}_A\,2 \cap \text{Ann}_A\,D_f$ using the standard basis of $A_f$. \qed

\section{$(\ast)$ is necessary}

\textbf{4.1 Lemma.} \textit{Let $x$ be in $A_f$, $X$ a subset of $S_n$ and $I$ an ideal of $A_f$ with $\sigma I = I$ for all $\sigma$ in $X$. We set $y := \prod_{\sigma\in X} \sigma x$. Then $I\cap \textnormal{Ann}_{A_f} x=0$ if and only if $I\cap \textnormal{Ann}_{A_f} y = 0.$}\\

\textbf{Proof.} Because $A_f \rightarrow A_f, y\mapsto \sigma y$ is an $A$-algebra homomorphism for all $\sigma$ in $S_n$, we may assume $1\in X$. Then $x$ divides $y$ and $\Leftarrow$ is trivial. We prove the implication $\Rightarrow$ by induction on $|X|$. For $|X| = 1$ we have $x=y$. For $|X|>1$ let $z\in I$ with $zy = 0$. Then
\begin{align*}
    0 = zy = (zx) \prod_{\sigma \in X\backslash{1}} \sigma x
\end{align*}
with $zx\in I$. By induction it follows that $zx = 0$ and thus $z = 0$. \qed\\

\newpage
\textbf{4.2 Lemma.} \textit{Let $1\leq i < j\leq n$ and set $\sigma := (i\, j)$.}
\begin{itemize}
    \item[(a)] \textit{Let $x\neq 0$ be a trivial $\sigma$-invariant. Then $x\tau_j \notin A_f^{\sigma,\textnormal{tr}}$.}
    \item[(b)] \textit{There is an inclusion $\textnormal{Ann}_{A_f} 2\cap \textnormal{Ann}_{A_f}(\tau_i+\tau_j)\subseteq A_f^{\sigma}$.}
\end{itemize}

\textbf{Proof.} First, let $i=1$ and $j=2$. We have $A_f = A[\tau_3, \tau_4, \dots, \tau_n]_{f_2}$. Thus, (b) follows by Lemma 3.1. Further, by Corollary 3.5
\begin{align*}
A_f^{\sigma, \text{tr}} = A_f^{S_2,\text{tr}} = A[\tau_3, \tau_4 \dots, \tau_n].
\end{align*}
Hence, (a) follows with the standard basis of $A_f$.\\
For arbitrary $i$ and $j$ we set $\pi := (1\,i)(2\,j)$. Then $\pi(\tau_1+\tau_2) = \tau_i+\tau_j$ and $\sigma = \pi (1\,2) \pi^{-1}$. Thus, (a) and (b) follow by the case $i=1$ and $j=1$ together with Remark 1.7 and Remark 2.5. \qed\\  

In the following we fix $\sigma:=(n-1\, n)$ and $G:=\langle S_{n-2}, \sigma \rangle$.\\

\textbf{4.3 Proposition.} \textit{If $\textnormal{Ann}_A 2\cap \textnormal{Ann}_A D_f \neq 0$, then there exists a non-trivial $\sigma$-invariant that is also a trivial $S_{n-2}$-invariant.}\\

\textbf{Proof.} By Lemma 4.1 it follows that $\text{Ann}_{A_f} 2 \cap \text{Ann}_{A_f} (\tau_{n-1} +\tau_n)\neq 0$. Using the standard basis of $A[\tau_n,\tau_{n-1}]_{f_{n-2}}=A_f$ we obtain
\begin{align*}
I:= \text{Ann}_{A[\tau_n,\tau_{n-1}]} 2 \cap \text{Ann}_{A[\tau_n,\tau_{n-1}]} (\tau_{n-1} +\tau_n) \neq 0.    
\end{align*}
By Corollary 3.5 we have $A_f^{S_{n-2},\text{tr}} = A[\tau_{n},\tau_{n-1}]_{f_{n-2}}^{S_{n-2},\text{tr}} = A[\tau_n,\tau_{n-1}]$ and thus $I \subseteq A_f^{S_{n-2},\text{tr}}$. By Lemma 4.2 (b) $I \subseteq A_f^\sigma$. Now let $0\neq x\in I$. If $x\notin A_f^{\sigma,\textnormal{tr}}$ we are done. Otherwise $\tau_n x\in I$ and $\tau_n x \notin A_f^{\sigma,\text{tr}}$ by Lemma 4.2 (a).\qed\\

\textbf{4.4 Lemma.} \textit{Let $x\in A_f^{S_{n-2},\textnormal{tr}}\cap A_f^\sigma$. Then $\Sigma_G^{S_n} x = x + y$ for $y \in A_f^{\sigma, \textnormal{tr}}$.}\\

\textbf{Proof.} Consider the (left) group action of $\sigma$ on $S_n/G$ and assume that the stabilizer of an element $\pi G$ is non-trivial. Then $\sigma \pi G = \pi G$ and so $\pi^{-1} \sigma \pi \in G$. Now write $\pi^{-1} \sigma \pi = (i\,j)$ with $\pi(i) = n-1$ and $\pi(j) = n$. It follows that $(i\,j)\in S_{n-2}$ or $(i\,j) = \sigma$. In the case $(i\,j) = \sigma$ we obtain $\pi \in G$. In the case $(i\,j)\in S_{n-2}$ we obtain $\pi x \in A_f^{\pi (i\,j) \pi^{-1},\text{tr}} = A_f^{\sigma,\text{tr}}$. In total this yields
\begin{align*}
    \Sigma_G^{S_n}x = x + \sum_{\pi} \pi x + \Sigma_1^\sigma z 
\end{align*}
for $z\in A_f$ where the sum runs over some $\pi \in S_n$ with $\pi x \in A_f^{\sigma,\text{tr}}$. Since also $\Sigma_1^\sigma z \in A_f^{\sigma,\text{tr}}$, we are done.\qed\\

The following corollary proves the necessity of $(\ast)$ in the main theorem.\\

\textbf{4.5 Corollary.} \textit{If $\textnormal{Ann}_A 2\cap \textnormal{Ann}_A D_f \neq 0$, then there exists a non-trivial \linebreak$S_n$-invariant.}\\

\textbf{Proof.} By Proposition 4.3 there exists a non-trivial $\sigma$-invariant $x$ that is also an $S_{n-2}$-invariant. Assume that $\Sigma_G^{S_n}x$ is a trivial $S_n$-invariant. By Lemma 4.4 $\Sigma_{G}^{S_n}x = x+y$ for a trivial $\sigma$-invariant $y$. This yields
\begin{align*}
x = \Sigma_G^{S_n}x -y\in A_f^{\sigma, \text{tr}}
\end{align*}
which is a contradiction. Hence, $\Sigma_G^{S_n}x$ is a non-trivial $S_n$-invariant. \qed

\section*{Acknowledgements}

I would like to express my gratitude to Professor Anne Henke who suggested this topic for my bachelor thesis, in which I developed the results of this paper. I am thankful to her for supervising this process.\\
I am also grateful to Professor Steffen Koenig for guidance while writing this article.

\small{University of Stuttgart, Institute of Algebra and Number Theory, Pfaffenwaldring 57, 70569 Stuttgart, Germany}\\
\small{\textit{E-mail address:} kevin.mobilmail@gmail.com}

\begin{thebibliography}{9}

\bibitem[EL]{EL} T. Ekedahl and D. Laksov, \emph{Splitting algebras, symmetric functions, and \hphantom GGalois theory}, J. Algebra Appl. 4 (2005), 59--75.

\bibitem[PZ]{PZ} M. Pohst and H. Zassenhaus, \emph{Algorithmic algebraic number theory, Ency-\hphantom Gclopedia of Mathematics and its Applications}, Cambridge University Press, \hphantom G1989.

\bibitem[Th]{Th} A. Thorup, \emph{On the invariants of the splitting algebra}, arXiv:1105.4478 \hphantom G[math.AC], 2011.

\end{thebibliography}
\end{document}